# Numerical solution of a certain type of integral equations on the real half-line

S. A. Belbas
Mathematics Department
University of Alabama
Tuscaloosa, AL 35487-0350. USA.

e-mail: SBELBAS@AS.UA.EDU

Abstracr. We develop a numerical method for solving a system of nonlinear integral equations involving two integral terms: at the current time t, one integral is taken from 0 to t, and a different integral is taken from t to infinity. We prove the convergence and the rate of convergence of our method. The discretization results in an infinite-dimensional nonlinear system, and we also prove results on the approximation of the solution of the infinite-dimensional system by solution of finite truncations.

## 1. Introduction and statement of the problem.

We are interested in the numerical solution of the following type of integral equation over the interval $0 \leq t < \infty$:

$$x(t) = x_0(t) + \int_0^t \exp(\alpha_1 s - \alpha_2 t) f(t,s,x(s)) \, ds + \int_t^\infty \exp(-\beta s + \gamma t) g(t,s,x(s)) \, ds$$

--- (1.1)

where the constants $\alpha_1, \alpha_2, \beta, \gamma$ satisfy

$$\alpha_1 \neq 0, \alpha_2 \geq 0, \alpha_1 \leq \alpha_2; \beta > 0, \beta \geq \gamma$$

--- (1.2)

and the functions f and g are bounded for an appropriate set of values of x and are Lipschitz in x; the exact properties of the functions f and g will be stated in the next section.

The unknown function x takes values in $\mathbb{R}^n$, and plainly so do the functions $x_0$, f and g.

Our motivation for studying (1.1) comes from integral equations over the real line ($-\infty < t < \infty$) of the form

$$z(t) = z_0(t) + \int_0^t F(t,s,z(s)) \, ds + \int_{-\infty}^t G(t,s,x(s)) \, ds$$

--- (1.3)

with two memory terms, one over an infinite interval and one starting at a finite time-instant (clearly, the exact value of 0 for the starting point in the first integral term is not crucial, and any finite time-instant as starting point could be transformed to an equation of the same form as (1.3)).

Now, (1.3) is numerically interesting for $t \leq 0$. Indeed, if a solution, say $z_1(\cdot)$, of (1.3) is known for all $t \leq 0$, then, for t>0, (1.3) becomes

$$z(t) = z_0(t) + \int_{-\infty}^0 G(t,s,z_1(s)) \, ds + \int_0^t \{F(t,s,z(s)) + G(t,s,z(s))\} \, ds$$

--- (1.4)

i.e. a standard Volterra integral equation for which numerical methods are well known (see, for example, [L]). For $t \leq 0$, by using the time-reversal transformation $\tau = -t$, $\sigma = -s$, and setting $Z(\tau) = z(-\tau)$, (1.3) takes the form

$$Z(\tau) = z_0(-\tau) - \int_0^\tau F(-\tau, -\sigma, Z(\sigma)) d\sigma + \int_\tau^\infty G(-\tau, -\sigma, Z(\sigma)) d\sigma$$

--- (1.5)

which is an equation of the form (1.1). Eq. (1.5) does not explicitly include the exponential weights that are shown in (1.1); the meaning of those exponential terms, when considered in the context of (1.5), is that they represent conditions on the rate of growth of the functions F and G of (1.5).

Particular cases of equations of the form (1.3) arise, in turn, from integro-differential systems of the form

$$\frac{du(t)}{dt} = f_1\left(t, u(t), \int_{-\infty}^t g_1(t, s, u(s)) ds\right); \; u(0) = u_0$$

--- (1.6)

Equations of the form of (1.6) appear in mathematical biology and in viscoelasticity theory.

By introducing a second unknown $v(t) := \int_{-\infty}^t g_1(t, s, u(s)) ds$, (1.6) becomes a system

$$u(t) = u_0 + \int_0^t f_1(s, u(s), v(s)) ds; \; v(t) = \int_{-\infty}^t g_1(t, s, u(s)) ds$$

--- (1.7)

which is a particular case of (1.3) with unknown function $[u^T \; v^T]^T$. It turns out that this particular case is numerically not easier than the full system (1.3) or, for that matter, not easier than (1.1). For this reason, we have examined the numerical solution of (1.1).

Usually, integral equations over infinite or semi-infinite intervals are numerically solved by a Nyström method that uses a Gauss-Laguerre or (depending on weights) other Gauss-related methods of numerical quadrature.(Nyström methods were introduced in [N], and a modern account of such methods may be found in [AH]; Gauss-Laguerre quadrature methods have been introduced in [S].) This is not feasible for (1.1), because of the simultaneous occurrence of two separate integrals over $[0,t]$ and over $[t,\infty)$. We get around this difficulty by using an infinite series approximation for the integral over $[t,\infty)$. This results in an infinite dimensional nonlinear system for the approximation of the function x(t) at the integration nodes. As a result, we are then faced with the problem of approximating the solution of this infinite system by using finite truncations; in section 3 of this paper, we solve the problem of convergence of the truncated solutions to the exact solution.

## 2. Discretization, and proof of convergence and rate of convergence.

We choose a discretization step h, and we are interested in approximating the solution x(t) of (1.1) at points $t_i := ih$, $i = 1, 2, \ldots$ . To this end, we approximate $\int_0^{ih} \exp(\alpha_1 s - \alpha_2 ih) f(ih, s, x(s)) \, ds$ by $\sum_{j=0}^{i-1} \left( \int_{jh}^{(j+1)h} \exp(\alpha_1 s - \alpha_2 ih) \, ds \right) f(ih, jh, x(jh))$ , and

$\int_{ih}^{+\infty} \exp(-\beta s + \gamma ih) g(ih, s, x(s)) \, ds$ by $\sum_{j=i}^{\infty} \left( \int_{jh}^{(j+1)h} \exp(-\beta s + \gamma ih) \, ds \right) g(ih, jh, x(jh))$.

We note the simple fact that

$$\int_{jh}^{(j+1)h} \exp(\alpha_1 s - \alpha_2 ih) \, ds = \frac{1}{\alpha_1} \exp((\alpha_1 j - \alpha_2 i)h)[\exp(\alpha_1 h) - 1],$$

$$\int_{jh}^{(j+1)h} \exp(-\beta s + \gamma ih) \, ds = \frac{1}{\beta} \exp((-\beta j + \gamma i)h)[1 - \exp(-\beta h)].$$

--- (2.1)

A consequence of (2.1) and the conditions on $\alpha_1$, $\alpha_2$, $\beta$, $\gamma$ stated in section 1 is that

$$\int_{jh}^{(j+1)h} \exp(\alpha_1 s - \alpha_2 ih) \, ds = O(h) \text{ uniformly in j for } 0 \leq j \leq i-1;$$

$$\int_{jh}^{(j+1)h} \exp(-\beta s + \gamma ih) \, ds = O(h) \text{ uniformly in i and j for } j \geq i .$$

--- (2.2)

The exact integral equation for the values x(ih) is

$$x(ih) = x_0(ih) + \sum_{j=0}^{i-1} \int_{jh}^{(j+1)h} \exp(\alpha_1 s - \alpha_2 ih) f(ih, s, x(s)) ds + \sum_{j=i}^{+\infty} \int_{jh}^{(j+1)h} \exp(-\beta s + \gamma ih) g(ih, s, x(s)) ds$$

--- (2.3)

We denote by $x_i$ the intended approximations to $x(ih)$. By utilizing the above approximate expressions for the integrals that appear in (2.3), we obtain the following equation for $x_i$, $i=1, 2, \ldots$ :

$$x_i = x_0(ih) + \sum_{j=0}^{i-1} \left( \int_{jh}^{(j+1)h} \exp(\alpha_1 s - \alpha_2 ih) ds \right) f(ih, jh, x_j) + \sum_{j=i}^{+\infty} \left( \int_{jh}^{(j+1)h} \exp(-\beta s + \gamma ih) ds \right) g(ih, jh, x_j)$$

--- (2.4)

In addition to the conditions stated in section 1, we assume the following:

(A1). The functions f and g are Lipschitz in x, for all x in a subset S of $\mathbb{R}^n$ to be specified below, uniformly in s and t, with Lipschitz constants $L_f$, $L_g$, respectively, and the Lipschitz constants satisfy $\frac{L_f}{|\alpha_1|} + \frac{L_g}{\beta} < 1$. Thus, for all s and t with $0 \leq s \leq t < +\infty$ and all x and y in S, we postulate

$|f(t,s,x) - f(t,s,y)| \leq L_f |x - y|$,

and for all s and t with $0 \leq t \leq s < +\infty$ and all x and y in S, we postulate

$|g(t,s,x) - g(t,s,y)| \leq L_g |x - y|$.

(A2). The functions f and g are continuous in all their arguments for x in S, and for $0 \leq s \leq t < +\infty$ for f, and respectively for $0 \leq t \leq s < +\infty$ for g, and for those domains of (s, t, x), we denote by $C_f$, $C_g$ the bounds for f and g, respectively:

$|f(t,s,x)| \leq C_f$, $|g(t,s,x)| \leq C_g$, for all (s, t, x) in the specified domains.

(A3). The set S is defined as the set of all n-dimensional vectors x that satisfy

$$|x| \leq \|x_0\|_\infty + \frac{C_f}{|\alpha_1|} + \frac{C_g}{\beta}.$$

(A4). The functions f and g are continuously differentiable in s and t, with bounded derivatives in the domains specified in (A2) above. Thus, we postulate that, for all (s, t, x) in the specified domains, we have

$$\left|\frac{\partial}{\partial t} f(t,s,x)\right| \leq D_f \;,\; \left|\frac{\partial}{\partial s} f(t,s,x)\right| \leq E_f \;,\; \left|\frac{\partial}{\partial t} g(t,s,x)\right| \leq D_g \;,\; \left|\frac{\partial}{\partial s} g(t,s,x)\right| \leq E_g$$

for some constants $D_f$, $D_g$, $E_f$, $E_g$.

(A5). The function $x_0(\cdot)$ is bounded and differentiable with bounded derivative, for all $t \geq 0$.

For every continuous bounded function $y(\cdot)$, we denote by Ty the function

$$(Ty)(t) := x_0(t) + \int_0^t \exp(\alpha_1 s - \alpha_2 ih) f(t,s,y(s)) ds + \int_t^\infty \exp(-\beta s + \gamma t) g(t,s,y(s)) ds \;;$$

--- (2.5)

for every bounded sequence $y = \{y_i : i = 1, 2, ...\} \in \ell_\infty(\mathbb{R}^n)$, we denote by $T_h y$ the sequence defined by

$$(T_h y)_i := x_0(ih) + \sum_{j=0}^{i-1} \left( \int_{jh}^{(j+1)h} \exp(\alpha_1 s - \alpha_2 ih) ds \right) f(ih, jh, y_j) +$$

$$+ \sum_{j=i}^{+\infty} \left( \int_{jh}^{(j+1)h} \exp(-\beta s + \gamma ih) ds \right) g(ih, jh, y_j)$$

--- (2.6)

(The norm in the space $\ell_\infty(\mathbb{R}^n)$ of all bounded $\mathbb{R}^n$-valued sequences is $\|y\|_\infty := \sup\{|y_i| : i = 1, 2, ...\}$.)

We have:

Lemma 2.1. The operator T maps continuous functions with values in S into continuous functions with values in S. The operator $T_h$ maps sequences with values in S into bounded sequences with values in S.

Proof: If $y(\cdot)$ is a continuous function with values in S, then we have

$$|(Ty)(t)| \leq \|x_0\|_\infty + \left(\int_0^t \exp(\alpha_1 s - \alpha_2 ih)\,ds\right) C_f + \sup_{s \geq t}\left(\int_t^\infty \exp(-\beta s + \gamma t)\,ds\right) C_g \leq$$

$$\leq \|x_0\|_\infty + \frac{C_f}{|\alpha_1|} + \frac{C_g}{\beta}$$

--- (2.7)

thus Ty takes values in S.
Next, if y is a sequence that takes values in S, then

$$|(T_h y)_i| \leq \|x_0\|_\infty + \left(\sum_{j=0}^{i-1} \int_{jh}^{(j+1)h} \exp(\alpha_1 s - \alpha_2 ih)\,ds\right) C_f +$$

$$+ \sup_{j \geq i}\left(\sum_{j=i}^\infty \int_{jh}^{(j+1)h} \exp(-\beta s + \gamma ih)\,ds\right) C_g =$$

$$= \|x_0\|_\infty + C_f \sum_{j=0}^{i-1}\left(\frac{1}{\alpha_1}\exp(\alpha_1 jh - \alpha_2 ih)[\exp(\alpha_1 h) - 1]\right) +$$

$$+ C_g \sup_{j \geq i}\sum_{j=i}^\infty \left(\frac{1}{\beta}\exp((-\beta j + \gamma i)h)[1 - \exp(-\beta h)]\right) =$$

$$= \|x_0\|_\infty + \frac{C_f}{|\alpha_1|}|\exp((\alpha_1 - \alpha_2)ih) - \exp(-\alpha_2 ih)| + \frac{C_g}{\beta}\exp(-(\beta - \gamma)ih) \leq \|x_0\|_\infty + \frac{C_f}{|\alpha_1|} + \frac{C_g}{\beta}$$

--- (2.8)

thus $T_h y$ also takes values in S. ///

Theorem 2.1. The operator T is a contraction (in the norm $\|\cdot\|_\infty$) on the set of continuous functions with values in S, and consequently (1.1) has a unique solution in that set of functions. The operator $T_h$ is a contraction (in the norm $\|\cdot\|_\infty$) on the set of sequences whose terms are in S, and consequently (2.4) has a unique solution in that set of sequences.

Proof: We have, for every two continuous functions y and z with values in S,

$$|(Ty)(t) - (Tz)(t)| \leq \left(\int_0^t \exp(\alpha_1 s - \alpha_2 ih)\,ds\right) L_f \|y - z\|_\infty$$

$$+ \sup_{s \geq t}\left(\int_t^\infty \exp(-\beta s + \gamma t)\,ds\right) L_g \|y - z\|_\infty \leq \left(\frac{L_f}{|\alpha_1|} + \frac{L_g}{\beta}\right)\|y - z\|_\infty = q\|y - z\|_\infty;$$

$$q := \frac{L_f}{|\alpha_1|} + \frac{L_g}{\beta} < 1$$

--- (2.9)

and a similar treatment, using the calculations contained in (2.8) above, applies to the operator $T_h$, which is thus established as contraction with the same value of the contraction constant q. ///

Now we state and prove a Gronwall-type inequality that will be needed for the next theorem. The exposition is facilitated by using an abstract formulation. In the proposition below, we avoid using the term "partially ordered Banach space" because the properties that we need for the proof are weaker than the properties that are specified in the standard definition of partially ordered Banach spaces. Although our formulation and proof seem natural, still we have not seen this formulation or proof of an abstract Gronwall-type inequality in the published literature. There exist several papers that prove a Gronwall-type inequality in a context that includes a partial order, but they do not contain any formulation or proof equivalent or close to ours. For this reason, we include a proof of this result.

Proposition 2.1. Let B be a Banach space over the field of real numbers with a partial order $\leq$ having the property that for all u, v, z, w in B, and for all nonnegative real numbers c, $\{u \leq v \text{ and } z \leq w\} \Rightarrow u + z \leq v + w$, and $u \leq v \Rightarrow cu \leq cv$, and also having the property that if two convergent (in norm) sequences $\{\mu_k : 1 \leq k < \infty\}$ and $\{v_k : 1 \leq k < \infty\}$ in B satisfy $\mu_k \leq v_k \ \forall k$, then $\lim_{k \to \infty} \mu_k \leq \lim_{k \to \infty} v_k$. Further, let Q be a (generally nonlinear) order-preserving operator on B and having the property that each iteration $Q^k$ is Lipschitz with Lipschitz constant $q_k$ with $\sum_{k=1}^\infty q_k < \infty$, and let $\delta$ be a vector in B that satisfies $\delta \leq Q\delta$. Let $\xi$ be the unique solution in B of the equation $\xi = Q\xi$. Then $\delta \leq \xi$.

Proof: Since Q has the property that each iteration $Q^k$ is Lipschitz with Lipschitz constant $q_k$ with $\sum_{k=1}^\infty q_k < \infty$, the equation $\xi = Q\xi$ can be solved by successive approximations starting with an arbitrary element of B, and we use $\delta$ as the initial vector, so that the successive approximations are $\xi_{(k)} = Q^k \delta$; on the other hand, the inequality $\delta \leq Q\delta$ and the order-preserving property of Q imply $Q\delta \le Q^2\delta$, thus $\delta \le Q\delta \le Q^2\delta$, and inductively $\delta \le Q\delta \le Q^2\delta \le \cdots \le Q^k\delta = \xi_{(k)}$; since $\xi_{(k)} \to \xi$ as $k \to \infty$, it follows that $\delta \le \xi$. ///

We remark that the condition that each $Q^k$ is Lipschitz with Lipschitz constant $q_k$ with $\sum_{k=1}^{\infty} q_k < \infty$, in the proposition above, includes as a particular case the case in which Q is a contraction.

Corollary 2.1. If Q is a linear contractive order-preserving operator on B, if $\delta$ satisfies $\delta \le C\lambda + Q\delta$ for some positive real constant C and some vector $\lambda$ in B, and if $\zeta$ solves $\zeta = \lambda + Q\zeta$, then $\delta \le C\zeta$.

Proof: $C\zeta$ satisfies $C\zeta = C\lambda + Q(C\zeta)$, and the conclusion follows by applying Proposition 2.1 to the mapping $\mu \mapsto C\lambda + Q\mu$ which is an order-preserving contraction on B. ///

Next, we establish the convergence and the rate of convergence of the solution of the approximate equation (2.4) to the exact solution of (1.1) as $h \to 0^+$.

Theorem 2.2. Let $\delta(h)$ be defined as

$$\delta(h) = \sup_{i \ge 0} |x(ih) - x_i|$$

where $x(\cdot)$ is the solution of (1.1) established in theorem 2.1, and $\{x_i : i = 1, 2, \ldots\}$ is the solution of (2.4) established in the same theorem. Then $\delta(h) = O(h)$.

Proof: For arbitrary but fixed step-size h, we set $\delta_i := |x(ih) - x_i|$. We note that, for each value of h, the sequence $\{\delta_i : i = 1, 2, \ldots\}$ is bounded since the two sequences $\{x(ih)\}$ and $\{x_i\}$ are bounded, and consequently $\delta(h)$ is finite. In order to simplify the notation, we also set

$$\varphi_j(h) := \int_{jh}^{(j+1)h} \exp(\alpha_1 s - \alpha_2 ih) ds, \quad \psi_{ij}(h) := \int_{jh}^{(j+1)h} \exp(-\beta s + \gamma ih) ds \quad (j \ge i).$$

As noted before, both $\varphi_i(h)$ and $\psi_{ij}(h)$ are of the order $O(h)$.

We have, from (2.3) and (2.4),

$$\delta_i \leq \sum_{j=0}^{i-1} \int_{jh}^{(j+1)h} \exp(\alpha_1 s - \alpha_2 ih) |f(ih, s, x(s)) - f(ih, jh, x_j)| ds +$$

$$+ \sum_{j=i}^{\infty} \int_{jh}^{(j+1)h} \exp(-\beta s + \gamma ih) |g(ih, s, x(s)) - g(ih, jh, x_j)| ds$$

--- (2.10)

We need an auxiliary result on the derivative of the exact solution of (1.1). By differentiating (1.1) with respect to t, we obtain

$$\dot{x}(t) = \dot{x}_0(t) + \exp((\alpha_1 - \alpha_2)t) f(t, t, x(t)) + \int_0^t \exp(\alpha_1 s - \alpha_2 t)[f_t(t, s, x(s)) - \alpha_2 f(t, s, x(s))] ds +$$

$$+ \int_t^{\infty} \exp(-\beta s + \gamma t)[\gamma g(t, s, x(s)) + g_t(t, s, x(s))] ds - \exp(-(\beta - \gamma)t) g(t, t, x(t))$$

--- (2.11)

According to the conditions (A1) through (A5), and for the unique solution x(t) of (1.1) that takes values in the set S, the terms
$\dot{x}_0(t), f(t, t, x(t)), g(t, t, x(t)), f_t(t, s, x(s)), g_t(t, s, x(s)), \exp(-\alpha t), \exp(-(\beta - \gamma)t)$

remain bounded, whereas the integrals $\int_0^t \exp(\alpha_1 s - \alpha_2 t) ds$, $\int_t^{\infty} \exp(-\beta s + \gamma t) ds$ are expressible in elementary ways and they are bounded for all nonnegative t. Thus the derivative $\dot{x}(t)$ remains bounded for all nonnegative t, and we set $C_0 := \sup_{t \geq 0} |\dot{x}(t)|$.

For the terms that appear on the right-hand side of (2.10), we have

$$|f(ih, s, x(s)) - f(ih, jh, x_j)| \leq |f(ih, s, x(s)) - f(ih, jh, x(s))| +$$
$$+ |f(ih, jh, x(s)) - f(ih, jh, x(jh))| + |f(ih, jh, x(jh)) - f(ih, jh, x_j)|.$$

--- (2.12)

In turn, the various terms on the right-hand side of (2.12) are estimated as follows:

$$|f(ih, s, x(s)) - f(ih, jh, x(s))| \leq ME_f(s - jh) \quad \text{for} \quad jh \leq s \leq (j+1)h$$

--- (2.13)

by the mean value inequality for vector-valued functions, where M is a constant that depends on the norm that is used in $\mathbb{R}^n$;

$$|f(ih, jh, x(s)) - f(ih, jh, x(jh))| \leq L_f |x(s) - x(jh)| \leq ML_f C_0 (s - jh) ;$$

--- (2.14)

$$|f(ih, jh, x(jh)) - f(ih, jh, x_j)| \leq L_f |x(jh) - x_j| = L_f \delta_j .$$

--- (2.15)

Of course, analogous inequalities hold for the function g. We set

$$C_1 := M(E_f + C_0 L_f), \quad C_2 := M(E_g + C_0 L_g).$$

Then we have, from (2.12) through (2.15),

$$\delta_i \leq C_1 \sum_{j=0}^{i-1} \int_{jh}^{(j+1)h} \exp(\alpha_1 s - \alpha_2 ih)(s - jh) ds + C_2 \sum_{j=i}^{\infty} \int_{jh}^{(j+1)h} \exp(-\beta s + \gamma ih)(s - jh) ds +$$

$$+ \sum_{j=1}^{i-1} L_f \left( \int_{jh}^{(j+1)h} \exp(\alpha_1 s - \alpha_2 ih) ds \right) \delta_j + \sum_{k=i}^{\infty} L_g \left( \int_{kh}^{(k+1)h} \exp(-\beta s + \gamma ih) ds \right) \delta_k$$

--- (2.16)

The summation in the next to last sum in (2.16) starts with j = 1 (instead of j = 0) because, by definition, $\delta_0 = 0$. We also use the standard convention that a sum over an empty set of indices is zero (to avoid treating the case i =1 as a special case).

By elementary calculations, it is verified that

$$C_1 \sum_{j=0}^{i-1} \int_{jh}^{(j+1)h} \exp(\alpha_1 s - \alpha_2 ih)(s - jh) ds + C_2 \sum_{j=i}^{\infty} \int_{jh}^{(j+1)h} \exp(-\beta s + \gamma ih)(s - jh) ds = O(h)$$
uniformly in i

--- (2.17)

and

$$\sum_{j=1}^{i-1} L_f \left( \int_{jh}^{(j+1)h} \exp(\alpha_1 s - \alpha_2 ih) ds \right) + \sum_{k=i}^{\infty} L_g \left( \int_{kh}^{(k+1)h} \exp(-\beta s + \gamma ih) ds \right) \leq \frac{L_f}{|\alpha_1|} + \frac{L_g}{\beta} = q < 1$$

--- (2.18)

Let $\{\zeta_i : i = 1,2,...\}$ be the unique bounded solution of

$$\zeta_i = 1 + \sum_{j=1}^{i-1} L_f \left( \int_{jh}^{(j+1)h} \exp(\alpha_1 s - \alpha_2 ih) ds \right) \zeta_j + \sum_{k=i}^{\infty} L_g \left( \int_{kh}^{(k+1)h} \exp(-\beta s + \gamma ih) ds \right) \zeta_k$$

--- (2.19)

The existence and uniqueness of a bounded solution of (2.19) follows from (2.18) and the contraction mapping fixed point theorem.

Then Corollary 2.1 implies that $\delta_i \leq O(h)\zeta_i = O(h)$ uniformly in $i$. This in turn implies the assertion of the theorem. ///

3. Numerical solution of the discretized integral equation.

Throughout this section, we keep the step-size h constant. We are interested in discovering ways to approximately solve the infinite-dimensional system, obtained in the previous section, for approximating the values x(ih) of the solution of the integral equation.

The system for the determination of $\{x_i : i = 1, 2, ...\}$ has the following form:

$$x_i = b_i + \sum_{j=0}^{i-1} \frac{1}{\alpha_1} \omega_{ij} F_{ij}(x_j) + \sum_{j=i}^{\infty} \frac{1}{\beta} \varphi_{ij} G_{ij}(x_j)$$

--- (3.1)

where the expressions $F_{ij}(x_j)$, $G_{ij}(x_j)$ are bounded by constants $C_f$, $C_g$, respectively, uniformly in i and j, and the terms $\omega_{ij}$, $\varphi_{ij}$ have the form

$$\omega_{ij} = \exp((\alpha_1 j - \alpha_2 i)h)[\exp(\alpha_1 h) - 1];$$
$$\varphi_{ij} = \exp((-\beta j + \gamma i)h)[1 - \exp(-\beta h)]$$

--- (3.2)

These quantities have (among others) the properties

$$\frac{1}{\alpha_1} \omega_{ij} \geq 0, \; \varphi_{ij} \geq 0, \; \sum_{j=0}^{i-1} \frac{1}{\alpha_1} \omega_{ij} \leq \frac{1}{|\alpha_1|}, \; \sum_{j=i}^{\infty} \varphi_{ij} \leq 1$$

--- (3.3)

The functions $F_{ij}$, $G_{ij}$ are Lipschitz in $x_j$ with Lipschitz constants $L_f$, $L_g$, respectively. One of our conditions is that $\frac{L_f}{|\alpha_1|} + \frac{L_g}{\beta} < 1$. Consequently, the system (3.1) can be solved, in principle, by a method of simple iterations in the space $\ell_\infty(\mathbb{R}^n)$. This is not directly useful for computer implementation, since actual implementation requires finite-dimensional operations. The obvious way would be to use a finite truncation of the infinite dimensional system. Such a truncation requires some analysis, since it is not obvious that the truncated solution would converge to the exact solution of (3.1) as the dimension of the truncated system goes to $\infty$. The convergence, and the topology in which convergence takes place, are not straightforward questions. One of the sources of the difficulty is the fact that the convergence of the series $\sum_{j=0}^{i-1} \frac{L_f}{\alpha_1} \omega_{ij} + \sum_{j=i}^{\infty} \frac{L_g}{\beta} \varphi_{ij}$ is not in general uniform in i; indeed, the first sum requires i terms to reach the value $\frac{L_f}{\alpha_1}[\exp((\alpha_1-\alpha_2)ih) - \exp(-\alpha_2 ih)]$, a value that is $\leq \frac{L_f}{|\alpha_1|}$, which is at a finite distance equal to $\frac{L_g}{\beta}\exp(-(\beta-\gamma)ih)$ from the exact value of the complete series, and, because we want to allow the possibility of $\beta = \gamma$, the last expression is in general no better than $\frac{L_g}{\beta}$. Thus, we distinguish two cases:

(I). $\beta > \gamma$.
(II). $\beta = \gamma$.

We denote by $\{x_{N,i} : i = 0, 1, 2, ..., N\}$ (with $x_{N,0} = x_0(0)$) the solution of the truncated system (corresponding to (5.1))

$$x_{N,i} = \sum_{j=0}^{(i-1)\wedge N} \frac{1}{\alpha_1} \omega_{ij} F_{ij}(x_{N,j}) + \sum_{j=i\wedge N}^{N} \frac{1}{\beta} \varphi_{ij} G_{ij}(x_{N,j})$$

--- (3.4)

The difference between the solution of the truncated system (5.4) and the complete system (3.1) satisfies

$$|x_{N,i} - x_i| \leq \sum_{j=0}^{(i-1)\wedge N} \frac{L_f}{\alpha_1} \omega_{ij} |x_{N,j} - x_j| + \sum_{j=i\wedge N}^{N} \frac{L_g}{\beta} \varphi_{ij} |x_{N,j} - x_j| + \sum_{j=N+1}^{\infty} \frac{C_f}{\beta} \varphi_{ij} |x_j|$$

$(i = 0, 1, 2, ..., N)$

--- (3.5)

We shall use the following notation: we denote by $(x_N)$ the extension of $x_N$ to an infinite-dimensional array by setting the entries after the N-th equal to zero, thus

$$(x_N)_i = \begin{cases} x_{N,i}, & \text{if } i \leq N \\ 0, & \text{if } i \geq N+1 \end{cases} ;$$

we also denote by $(x)_N$ the truncation of the infinite-dimensional array $x$ to its first N entries, thus

$$(x)_{N,i} = \begin{cases} x_i, & \text{if } i \leq N \\ 0, & \text{if } i \geq N+1 \end{cases} .$$

For case (I), we have:

Theorem 3.1. If $\beta > \gamma$, then the solution of the truncated system (3.4) converges to the exact solution x of (3.1) in the sense that

$$\lim_{N \to \infty} \| x_N - (x)_N \|_\infty = 0 .$$

Proof: For every $\varepsilon > 0$, we choose $N_1 = N_1(\varepsilon)$ so that $\dfrac{L_g}{\beta} \exp(-(\beta - \gamma)ih) < \varepsilon$ whenever $i > N_1$. (This is possible because $\exp(-(\beta - \gamma)ih) \to 0$ as $i \to \infty$.) Next, we note that, for $i \leq N$, we have

$$\sum_{j=N+1}^{\infty} \exp(-\beta jh + \gamma ih)[1 - \exp(-\beta h)] = \exp(-\beta(N+1)h + i\gamma h) \to 0 \text{ as } N \to \infty .$$

Consequently, we can choose an $N_2 = N_2(\varepsilon) > N_1(\varepsilon)$ so that, for all i that satisfy $0 \leq i \leq N_1$, if $N \geq N_2(\varepsilon)$, then $\sum_{j=N+1}^{\infty} \dfrac{C_g}{\beta} \varphi_{ij} < \varepsilon$. By truncating the infinite dimensional system (3.1) to terms with $0 \leq i \leq N_2$, we find, by dint of (3.5),

$$|x_{N_2,i} - x_i| \leq \sum_{j=0}^{(i-1)\wedge N_2} \dfrac{L_f}{\alpha_1} \omega_{ij} |x_{N_2,j} - x_j| + \sum_{j=i \wedge N_2}^{N_2} \dfrac{L_g}{\beta} \varphi_{ij} |x_{N_2,j} - x_j| + \sum_{j=N_2+1}^{\infty} \dfrac{C_g}{\beta} \varphi_{ij} |x_j|$$

$(i = 0, 1, 2, ..., N_2)$

--- (3.6)

and therefore

$$\| x_{N_2} - (x)_{N_2} \|_\infty \leq q \| x_{N_2} - (x)_{N_2} \|_\infty + \varepsilon \| x \|_\infty$$

from which

$$\| x_{N_2} - (x)_{N_2} \|_\infty \leq \dfrac{\varepsilon}{1-q} \| x \|_\infty$$

--- (3.7)

and (3.7) is tantamount to the assertion of the theorem. ///

At this point, it is useful to prove the following:

Proposition 3.1. If $\beta > \gamma$ and $\alpha_1 < \alpha_2$, if $\lim_{i \to \infty} b_i = 0$, then also $\lim_{i \to \infty} x_i = 0$.

Proof: We obtain, from (3.1),

$$|x_i| \leq |b_i| + \frac{C_f}{\alpha_1} \sum_{j=0}^{i-1} \omega_{ij} + \frac{C_g}{\beta} \sum_{j=i}^{\infty} \varphi_{ij} = |b_i| + \frac{C_f}{\alpha_1}[\exp((\alpha_1 - \alpha_2)ih) - \exp(-\alpha_2 ih)] + \frac{C_g}{\beta}\exp(-(\beta-\gamma)ih)$$

--- (3.8)

and every term on the right-hand side of (3.8) goes to zero as i goes to infinity. ///

In the case $\beta = \gamma$, we use a transformation $y_i = e^{-i\delta h} x_i$, $i = 0,1,2,...$ where $\delta$ is a positive number, that satisfies certain conditions to be specified below. The Lipschitz condition for the new variables $y_i$ is obtained by substituting $\exp(j\delta)y_j$ in lieu of $x_j$ in $F_{ij}$ and $G_{ij}$:

$$|F_{ij}(\exp(i\delta h)z_j) - F_{ij}(\exp(i\delta h)w_j)| \leq L_f \exp(j\delta h)|z_j - w_j|;$$
$$|G_{ij}(\exp(i\delta h)z_j) - G_{ij}(\exp(i\delta h)w_j)| \leq L_g \exp(j\delta h)|z_j - w_j|$$

--- (3.9)

The infinite dimensional vector y solves

$$y_i = \exp(-i\delta h)\{b_i + \sum_{j=0}^{i-1} \frac{\omega_{ij}}{\alpha_1} F_{ij}(\exp(j\delta h)y_j) + \sum_{j=i}^{\infty} \frac{\varphi_{ij}}{\beta} G_{ij}(\exp(j\delta h)y_j)\}$$

--- (3.10)

If we denote by $U_i(y)$ the right-hand side of (3.10), then $U_i$ is Lipschitz in y with Lipschitz constant

$$\frac{L_f}{\alpha_1} \sum_{j=0}^{i-1} \omega_{ij} \exp((j-i)\delta h) + \frac{L_g}{\beta} \sum_{j=i}^{\infty} \varphi_{ij} \exp((j-i)\delta h) \ .$$

We have

$$\omega_{ij} \exp((j-i)\delta h) = \exp((\alpha_1 + \delta)jh - (\alpha_2 + \delta)ih)[\exp(\alpha_1 h) - 1]$$

and, by choosing $\delta$ so small that $\alpha_1 + \delta$ has the same sign as $\alpha_1$, we obtain the quantities $\alpha_1' := \alpha_1 + \delta$, $\alpha_2' := \alpha_2 + \delta$ which satisfy the same conditions as those postulated for the original quantities $\alpha_1, \alpha_2$. We note that

$$\sum_{j=0}^{i-1} \frac{1}{\alpha_1} \omega_{ij} \exp((j-i)\delta h) = \frac{1}{\alpha_1}[\exp((\alpha_1-\alpha_2)ih) - \exp(-(\alpha_2+\delta)ih)]\frac{\exp(\alpha_1 h)-1}{\exp((\alpha_1+\delta)h)-1} \leq$$

$$\leq \frac{1}{|\alpha_1|} \frac{\exp(\alpha_1 h)-1}{\exp((\alpha_1+\delta)h)-1} .$$

--- (3.11)

We also have

$$\varphi_{ij} \exp((j-i)\delta h) = \exp(-(\beta-\delta)jh + (\gamma-\delta)ih)[1-\exp(-\beta h)]$$

and now we choose $\delta$ so small that all previous properties hold and additionally $\beta-\delta$, $\gamma-\delta$ remain positive. We note that

$$\sum_{j=i}^{\infty} \varphi_{ij} \exp((j-i)\delta h) = \exp(-(\beta-\gamma)ih)\frac{1-\exp(-\beta h)}{1-\exp(-(\beta-\delta)h)} \leq$$

$$\leq \frac{1-\exp(-\beta h)}{1-\exp(-(\beta-\delta)h)} .$$

--- (3.12)

The expressions $\dfrac{\exp(\alpha_1 h)-1}{\exp((\alpha_1+\delta)h)-1}$ and $\dfrac{1-\exp(-\beta h)}{1-\exp(-(\beta-\delta)h)}$ that appear on the right-hand sides of (3.11) and (3.12) are both independent of i and go to 1 as $\delta$ goes to 0. Therefore, it is possible to choose $\delta>0$ so small that all the previous requirements are satisfied and in addition

$$\frac{L_f}{|\alpha_1|} \frac{\exp(\alpha_1 h)-1}{\exp((\alpha_1+\delta)h)-1} + \frac{L_g}{\beta} \frac{1-\exp(-\beta h)}{1-\exp(-(\beta-\delta)h)} < 1 .$$

--- (3.13)

We set $\vartheta := \dfrac{L_f}{|\alpha_1|} \dfrac{\exp(\alpha_1 h)-1}{\exp((\alpha_1+\delta)h)-1} + \dfrac{L_g}{\beta} \dfrac{1-\exp(-\beta h)}{1-\exp(-(\beta-\delta)h)}$. For simplicity, we represent the system satisfied by y in the following form:

$$y_i = b_i(\delta) + \frac{1}{\alpha_1} \sum_{j=0}^{i-1} \omega_{ij}(\delta) \Phi_{ij}(y_j) + \frac{1}{\beta} \sum_{j=i}^{\infty} \varphi_{ij}(\delta) \Psi_{ij}(y_j)$$

--- (3.14)

where

$b_i(\delta) := \exp(-i\delta h) b_i$, $\omega_{ij}(\delta) := \omega_{ij} \exp((j-i)\delta h)$, $\varphi_{ij}(\delta) := \varphi_{ij} \exp((j-i)\delta h)$,

$\Phi_{ij}(y_j) := F_{ij}(\exp(j\delta h) y_j)$, $\Psi_{ij}(y_j) := G_{ij}(\exp(i\delta h) y_j)$ .

The operator appearing on the right-hand side of (3.14) is a contraction on $\ell_\infty(\mathbb{R}^n)$ with contraction constant $\vartheta<1$ defined above. Thus (3.14) has all the essential properties of the original system (3.1) and the additional property $\lim_{i\to\infty} y_i = 0$ (since $y_i = e^{-i\delta h} x_i$ and $x_i$ is bounded uniformly in i). We shall use the same notation for truncated systems corresponding to (3.14) as for the truncations of (3.1).

We have:

Theorem 3.2. We assume that $\delta$ is sufficiently small so that all the properties detailed above are satisfied. Then the truncated solutions $y_N$ corresponding to the system (3.14) converge uniformly to the exact solution y, in the sense that the extensions $(y_N)$ satisfy $\lim_{N\to\infty} \|(y_N) - y\|_\infty = 0$.

Proof: We have (by subtracting the equations satisfied by y and $(y_N)$, in a way similar to deriving (3.5))

$$|(y_N)_i - y_i| \leq \sum_{j=0}^{(i-1)\wedge N} \frac{L_f}{\alpha_1} \omega_{ij}(\delta) |(y_N)_i - y_i| + \sum_{j=i\wedge N}^{N} \frac{L_g}{\beta} \varphi_{ij}(\delta) |(y_N)_i - y_i| + \sum_{j=i\vee(N+1)}^{\infty} \frac{L_f}{\beta} \varphi_{ij}(\delta) |y_j|.$$

--- (3.15)

The quantity $\sum_{j=i\vee(N+1)}^{\infty} \varphi_{ij}(\delta)$ is bounded uniformly in i, because $\sum_{j=i\vee(N+1)}^{\infty} \varphi_{ij}(\delta) \leq \sum_{j=i}^{\infty} \varphi_{ij}(\delta)$ and $\sum_{j=i}^{\infty} \varphi_{ij}(\delta)$ is bounded uniformly in i by (3.12). Since also $y_j \to 0$ as $j\to\infty$, we conclude that $\lim_{N\to\infty} \sum_{j=i\vee(N+1)}^{\infty} \frac{L_f}{\beta} \varphi_{ij}(\delta) |y_j| = 0$. Also, (3.14) implies

$$\|(y_N) - y\|_\infty \leq \vartheta \|(y_N) - y\|_\infty + \sum_{j=i\vee(N+1)}^{\infty} \frac{L_f}{\beta} \varphi_{ij}(\delta) |y_j|$$

and consequently $\|(y_N) - y\|_\infty \leq \frac{1}{1-\vartheta} \sum_{j=i\vee(N+1)}^{\infty} \frac{L_f}{\beta} \varphi_{ij}(\delta) |y_j| \to 0$ as $N\to\infty$ . ///

Corollary 3.1. Under the conditions of theorem 3.2, for every $\varepsilon > 0$ and for every i there is an $N(\varepsilon,i)$ such that $|x_{N,i} - x_i| < \varepsilon$ whenever $N \geq N(\varepsilon,i)$.

Proof: For arbitrary but fixed i, we choose $N(\varepsilon,i) \geq i$ and also so large that, whenever $N \geq N(\varepsilon,i)$, we have $|y_{N,i} - y_i| < \exp(-i\delta h)\varepsilon$. This is possible because, as a consequence of theorem 3.2, $y_{N,i} \to y_i$ uniformly in i, as $N \to \infty$. Since $x_i = \exp(i\delta h)y_i$ and $x_{N,i} = \exp(i\delta h)y_{N,i}$, it follows that, for $N \geq N(\varepsilon,i)$, we have $|x_{N,i} - x_i| < \varepsilon$. ///